
\magnification=\magstep1
\documentstyle{amsppt}

\def\pr{\noindent{\sl Proof. }}
\def\re{\noindent{\sl Remark. }}
\def\E{{\Cal E}} 
\def\F{{\Cal F}}
\def\G{{\Cal G}}

\def\I{{\Cal I}}
\def\J{{\Cal J}}

\def\t{\tilde}
\def\O{{\Cal O}}
\def\P{{\Bbb P}}
\def\Z{{\Bbb Z}}
\def\Tors{{\operatorname {Tors}}}
\def\Ext{{\operatorname{Ext}}}

\def\mult{{\operatorname{mult}}}

\def\ra{\longrightarrow}
\def\raa{\rightarrow}

\def\ot{\otimes}
\def\op{\oplus}

\def\today{\ifcase\month \or January \or February \or March \or April
\or May \or June \or July \or August \or September \or October
\or November \or December \fi \space \number\day,  \number \year}

\def\:{\colon}

\topmatter
\title 
A note on $k$-jet ampleness on surfaces
\endtitle
\author Adrian Langer
\endauthor

\rightheadtext{jet ampleness}
{\date{\today}\enddate}

\address{ Adrian Langer:
Instytut Matematyki UW, 
ul.~Banacha 2, 02--097 Warszawa, Poland}
\endaddress
\email{alan\@mimuw.edu.pl}\endemail

\abstract We prove Reider type criterions for $k$-jet spannedness
and $k$-jet ampleness of adjoint bundles for surfaces with
at most rational singularities. Moreover, we prove that on smooth
surfaces $[n(n+4)/4]$-very ampleness implies $n$-jet ampleness. 
\endabstract

\thanks The author was partially supported by Foundation 
for Polish Science
\endthanks


\endtopmatter

\document 
\vskip .3 cm

\heading  Introduction\endheading

Let $L$ be a Cartier divisor on a normal projective surface $X$ and 
$k$ a non-negative integer.
$L$ generates $k$-jets at a point $x\in X$ if the restriction map 
$H^0(L)\to H^0(L\ot \O_X/m_x^{k+1})$ is onto.
$L$ is $k$-jet generated (or $k$-jet spanned) if it generates
$k$-jets at each point of $X$.
$L$ is $k$-jet ample if for any distinct points $x_1$,..., $x_r$ in $X$
and positive integers $k_1$,..., $k_r$ with $k_1+...+k_r=k+1$
the restriction map 
$H^0(L)\to H^0(L\ot \O_X/m_{x_1}^{k_1}\ot ...\ot m_{x_r}^{k_r})$ 
is onto.
Note that $L$ is $0$-jet ample if and only if it is $0$-jet generated and
if and only if it is spanned by global sections.

The main aim of this note is to establish Reider type criterions 
for $k$-jet spannedness and $k$-jet ampleness
of an adjoint bundle. Up to now there were a few trials 
to find such criterions
(see, e.g., [BS2] and [Laz], Section 7) but the optimal results
were not known.
A Reider type criterion is well known in the case
of $k$-very ampleness (see [BS1, Theorem 2.1]) and it implies
a weak form of the criterion for $k$-jet ampleness
(see [BS2, Proposition 2.1]). 
A different version of the criterion 
for $k$-jet spannedness with better bounds
was proved in [Laz, Theorem 7.4].

The paper is divided into 4 sections. 
In the first section we recall some results used in the paper.
In Section 2 we prove that on smooth surfaces already $[{n(n+4)\over 4}]$-very
ampleness implies $n$-jet ampleness. This together with [BS1, Theorem 2.1]
gives much better Reider type criterion for $k$-jet ampleness 
than those mentioned above.
In Section 3 we give a direct proof of even better criterion for $k$-jet
spannedness and $k$-jet ampleness (without using results of [BS1])
describing also the boundary case in terms of the Seshadri constant. 
This version of the Reider type criterion is new
even in the case of $0$-jet spannedness (i.e., in the classically known
globally generated case). Our proof works also in a larger category
of normal surfaces with at most rational singularities whereas the 
results of [BS1], [BS2] and [Laz] are known only for smooth surfaces.
In particular, our theorem holds for canonical surfaces, where it was
used to prove that $|2K_X|$ has no base components for surfaces
with $K_X^2=4$, $p_g=q=0$ (see [La2], Theorem 0.1). 
In the last section we try to explain (after [Laz]) how Seshadri constants
appear in the study of adjoint bundles and we give an example of 
$1$-jet spanned but not $1$-jet ample line bundle.

\heading 1. Preliminaries \endheading

{\bf 1.1.}
Let $L$ be a line bundle on a normal surface $X$. It is very natural
to consider the following definition: $L$ is called 
{\sl $k$-point generated} 
(or {\sl $k$-point spanned}) at a point $x$ if the restriction 
map $H^0(\O_X(L))\to H^0(\O_\zeta(L))$ is onto for any cluster $\zeta$ 
supported on $x$ and of degree $\le k+1$. In fact it would be more
natural to call it $k$-generated, but this notion is reserved for 
something slightly different.
Recall also that $L$ is called {\sl $k$-very ample} if the restriction map
$H^0(\O_X(L))\to H^0(\O_\zeta(L))$ is onto for any degree $\le k+1$
cluster $\zeta$ in $X$.

We will use those notions in Section 2.

\medskip

\proclaim{Lemma 1.2}
Let $D$ be a Weil divisor on a normal surface. If $D^2\ge 0$ and $DL>0$
for some nef divisor $L$, then $D$ is pseudoeffective. Moreover, 
$D$ is big unless $D^2=0$ and $D$ is nef.
\endproclaim

\noindent
{\sl Sketch of the proof.}
The first part of the lemma follows easily from the Hodge index theorem.
The second one follows from the Zariski decomposition for $D$, Q.E.D.

\medskip
The following lemma is a slightly modified version of 
Corollary 3.7, [La1], and can be obtained similarly as in [La1]
by using Lemma 1.2.

\proclaim{Lemma 1.3} 
Let $X$  be a normal projective surface and $L$ a pseudoeffective
Weil divisor on $X$. If $L^2>0$ or $L^2=0$ and $L$ is not nef
then every nontrivial extension $\E\in \Ext ^1(\O(K_X+L),\omega_X)$
is Bogomolov unstable.
\endproclaim

Let us also recall the definition of the Seshadri constant.

\proclaim{Definition 1.4}
Let $X$ be a normal surface with at most rational singularity at a  point $x$
and let $f\:Y\to X$ be the minimal resolution of the singularity at $x$
(or   the blow up of $X$ at $x$ if $x$ is smooth).
Let $Z$ denote the fundamental cycle (respectively: the exceptional divisor).
A Seshadri constant of a divisor $L$ at  $x$ is defined  as
$$\epsilon (L,x)=\sup\{\epsilon\ge 0| f^*L-\epsilon\cdot Z{\hbox{ is nef }}
\},$$
whenever it is a well defined real number.
\endproclaim

\heading 2. Relations with $k$-very ampleness \endheading

Let $x$ be a smooth point of a surface $X$.
Our first problem will be to determine the number 
$$l_n=\max
\{\deg\zeta\: \, \zeta{\hbox{  is Gorenstein and }m_x^{n+1}\subset\I_\zeta}
\hbox{ in } \O_{X,x} \} .$$

Clearly,  $l_0=1$ and $l_n< {n+2 \choose 2}$ for $n\ge 1$ since
$m_x^{n+1}$ is 
not locally generated by two elements and hence the corresponding
cluster is not Gorenstein. In fact, we have the following theorem:

\proclaim{Theorem 2.1}
$$l_n=\left[ {(n+2)^2\over 4}\right]$$
\endproclaim

\pr
Note that $l_n\ge [ {(n+2)^2\over 4}]$ since if $z_1$, $z_2$ are local 
parameters and $\I_\zeta= (z_1^{k+1},z_2^{k+1})$ for $n=2k$, or
$\I_\zeta= (z_1^k,z_2^{k+1})$ for $n=2k-1$, then $\zeta$ is Gorenstein
and $\deg \zeta=[{(n+2)^2\over 4}]$.

Let $(X,x)$ be a germ of a smooth surface. By the above it is sufficient 
to prove that if for a cluster $\zeta$ supported on the point $x$,
the ideal $\I_\zeta$ is generated by two elements $f$ and $g$
and it contains $m_x^{n+1}$ then  $\deg\zeta\le {(n+2)^2\over 4}$
(let us recall that in the smooth case Gorenstein cluster on
a surface means locally complete intersection).

Let $V_1$ and $V_2$ be curves defined by $f$ and $g$, respectively.
It is easy to see that there exists a sequence of blow ups
$\pi \:X_k{\mathop{\ra}^{p_k}}X_{k-1}\ra ...\ra X_1{\mathop
{\ra}^{p_1}}X$
such that 
$$\pi ^*V_i={\t V_i}+\sum e_{Q_j}(V_i) E_j,$$
where $\t V_i$ is a strict transform of $V_i$, $E_j$
denotes (a pull back of) an exceptional divisor of $p_j$, 
$Q_j=p_j(E_j)$, $e_{Q_j}(V_i)$ denotes a multiplicity of $V_i$ along $Q_j$,
$\t V_1$ and $\t V_2$ are disjoint and 
the divisor $\t V_1+\t V_2+\sum E_j$ is normal crossing. 

Now note that the scheme-theoretical preimage of $I_\zeta$, i.e.,
$(\pi^*f,\pi^*g)=\O(-\t V_1-\sum e_{Q_j}(V_1)E_j)+\O(-\t V_2-
\sum e_{Q_j}(V_2)E_j)$
does not contain $\J=\O(-\sum(e_{Q_j}(V_1)+e_{Q_j}(V_2))E_j-2E_k)$
and $\pi_*\J\subset m_x^{\sum(e_{Q_j}(V_1)+e_{Q_j}(V_2))-2}$.

Therefore if $m_x^{n+1}\subset (f,g)$ then $n\ge \sum(e_{Q_j}(V_1)
+e_{Q_j}(V_2))-2$.

By the definition the degree of $\zeta$ is a local 
intersection multiplicity of $V_1$ and $V_2$ and from the above it is easy 
to see that
$$\deg\zeta=\sum e_Q(V_1)e_Q(V_2),$$ 
where we  sum  over all infinitely near points 
$Q$ of $X$
(this is the formula of M. Noether; see [Fu], Example 12.4.2). 
Now the theorem follows from the obvious inequality
$$\sum e_Q(V_1)e_Q(V_2) \le {\left(\sum (e_Q(V_1)+e_Q(V_2)\right)^2\over 4}
\le {(n+2)^2\over 4}.
$$

\medskip

\proclaim{Theorem 2.2}
Let $L$ be a Cartier divisor on a normal surface $X$ and $x$ a smooth
point of $X$.
\item{1.} If $L$ is $(l_k-1)$-point generated at a point $x$ then 
$L$ is $k$-jet generated at $x$.
\item{2.} If $X$ is smooth and $L$ is $(l_k-1)$-very ample 
then $L$ is $k$-jet ample.
\endproclaim

\pr
By the cohomology exact sequence 
$$
H^0(\O (L))\ra H^0(\O (L)/m_x^{k+1})\ra H^1(m_x^{k+1}\O (L))\ra H^1(\O (L))
$$
and $\Ext ^1(m_x^{k+1}\O (L),\omega _X)$ ($\Ext ^1(\O (L),\omega _X)$)
is dual to $H^1(m_x^{k+1}\O (L))$ ($H^1(\O (L))$, respectively)
by the Serre duality theorem. 
If $L$ is not $k$-jet generated  at $x$ then using this  
one can see that there exists an extension 
$\E\in \Ext ^1(m_x^{k+1}\O (L),\omega _X)$ not coming
from $\Ext ^1(\O(L),\omega _X)$. Set $\F=\E^{**}$.
Then $\F\in \Ext ^1(\I_\zeta \O(L),\omega _X)$
for some cluster $\zeta$ contained in the cluster $\O_X/m_x^{k+1}$.

Since $\F$ is locally free, $\zeta$ is a locally complete intersection
and hence $\deg \zeta\le l_k$. But $L$ is $(l_k-1)$-point generated at $x$,
so by the similar arguments as above we prove that every extension in 
$\Ext ^1(\I_\zeta \O(L),\omega _X)$ comes from
$\Ext ^1( \O(L),\omega _X)$, a contradiction.

The proof in the ample case is very similar once we know that
$$l_{k_1-1}+...+l_{k_r-1}\le l_k\leqno (*)$$
for any positive integers $k_1$,..., $k_r$ with $k_1+...+k_r=k+1$.
But we have
$$\sum (k_i+1)^2=\sum k_i^2 +2(k+1)+r\le \sum k_i^2+\sum _{i\ne j}k_i k_j+1+
2(k+1)=$$
$$=(k_1+...+k_r)^2+2(k+1)+1=(k+2)^2$$
from which the inequality $(*)$ follows.
This finishes the proof of the theorem.

\medskip

\noindent
{\sl Remarks.}

\noindent
(1) Similar in vein but weaker theorem  was proved in any dimension by
Beltrametti and Sommese (see [BS2, Proposition 2.1]). It would be interesting
to know whether our theorem can be also  generalized to the higher dimensional
case.

\noindent 
(2) Note that already the trivial bound $l_n< {n+1 \choose 2}$ for $n>1$ 
implies in the surface case better theorem than Proposition 2.1, [BS2]. 
In fact, Theorems 2.1 and 2.2 imply a part of Corollary 3.2 
(Theorem 3.4) by using Reider type theorem 
for $k$-point spannedness ($k$-very ampleness, respectively; see [BS1]).

\heading 3. Main theorems \endheading

\proclaim{Theorem 3.1}
Let $L$ be a pseudoeffective Weil divisor on a normal 
projective surface $X$.
Assume that $K_X+L$ is Cartier and generates $(k-1)$-jets but not $k$-jets
at a smooth point $x$.
\item{(1)} If $L^2>(k+2)^2$ then there exists a curve $D$ containing
a Gorenstein cluster $\zeta$ such that $\I_\zeta$ contains $m_x^{k+1}$
but it does not contain  $m_x^k$ and such that the map 
$H^0(\O_D(K_X+L))\to H^0(\I _\zeta\O_D(K_X+L))$ is not onto.
In particular, $|\O_D(K_X+L)|$ does not  generate $k$-jets at the point $x$.
Moreover, $L-2D$ is pseudoeffective, numerically nontrivial and the following
inequalities hold:
$$LD-\deg \zeta\le 2p_aD-2-K_XD$$
and
$$LD-{1\over 4}(k+2)^2\le D^2.$$ 
\item{(2)} If $L^2=(k+2)^2$, then either there exists a curve $D$
as in (1) or $\epsilon (x, L)=k+2$.
\endproclaim

\pr
Since $K_X+L$ does not generate $k$-jets at $x$ and generates $(k-1)$-jets,
by the Serre duality theorem there exists a nontrivial extension 
$\F\in \Ext ^1(m_x^{k+1}\O(K_X+L), \omega _X)$ not lying in the image
of the natural map $\Ext ^1(m_x^{k}\O(K_X+L), \omega _X)\to
\Ext ^1(m_x^{k+1}\O(K_X+L), \omega _X)$.
Let $\E$ be a reflexivisation of $\F$. Then the cokernel of the natural map
$\omega_X\to \E$ twisted by $\O(-K_X-L)$ defines an ideal of a cluster 
$\zeta$. Clearly, $m^{k+1}_x\subset \I_\zeta$ but  $m_x^k \not\subset I_\zeta$.
Since $\E\in \Ext ^1(\I_\zeta\O(K_X+L), \omega _X)$ is reflexive,
$\zeta$ is a Gorenstein scheme.

Let $f\: Y\to X$ be the blow up of $X$ at $x$ and let us denote by $E$
the exceptional divisor. Then $f_*\O(-(k+1)E)=m_x^{k+1}$ and therefore
there exists a rank $2$ reflexive sheaf
$\F'\in \Ext^1(\O_Y(f^*(K_X+L)-(k+1)E),\omega_Y)$ such that 
$f_*\F'=\F$. 
By Lemma 1.2 and the assumption that $L^2\ge (k+2)^2$,
the sheaf $\F'$ is Bogomolov unstable unless $f^*L-(k+2)E$ is nef
and $L^2=(k+2)^2$ in which case $\epsilon (L,x)=k+2$.
Therefore $\E=(f_*\F')^{**}$ is also Bogomolov unstable reflexive sheaf.
Let $\O(K_X+A)$ be a maximal destabilizing subsheaf of $\E$.
Then the composition map  $\O(K_X+A)\to \E\to \I_\zeta\O(K_X+L)$
is nonzero and its image twisted by $\O(-K_X-L)$ defines an ideal
of effective divisor $D$ containing $\zeta$ and linearly equivalent to 
$L-A$. By the instability of $\E$ the divisor $L-2D$ 
is pseudoeffective and numerically nontrivial.

Using rather lengthy arguments with diagram chasing
one can prove that $\zeta$ is in a (very) special position 
with respect to $\O_D(K_X+L)$; see Lemma 4.4.1, [La3], or 
the proof of Theorem 4.7, [La1]. In particular, it follows that
there is an injection $\I_\zeta\O _D(K_X+L)\hookrightarrow \omega_D$ and hence
$2p_aD-2\ge  D(K_X+L)-\deg \zeta$ which is equivalent to
$DL-\deg\zeta\le D^2$.
It also shows that $\O_D(K_X+L)$ is not $k$-jet generated
at $x$. On the other hand it is $(k-1)$-jet generated at $x$, since
$K_X+L$ is $(k-1)$-jet generated at $x$.

The other inequality can be proven in much the same way as
a similar inequality in Theorem 4.7, [La1], Q.E.D.

\medskip
\re
There are some variants of the theorem which seems to be worth of pointing
out:
\item{1.} If under the assumptions of Theorem 3.1 the surface $X$
is  Gorenstein one can get also that $LD-\deg\zeta\le D^2$
which is better than the second inequality of Theorem 3.1
(and it is not equivalent to the first one unless $X$ is smooth).
This inequality  is equivalent to  $c_2(\E)\ge (K_X+A)(K_X+D)$ 
which follows from an exact sequence
$0\raa\O (K_X+A)\ra\E \ra\G\raa 0,$
where $\G^{**}=\O(K_X+D)$
(note  that this is more difficult  than it seems since $K_X+A$ and $K_X+D$
are not necessarily Cartier; nevertheless it is still true; see, e.g., 
[La3], Proposition 2.15).
\item{2.} If $X$ is smooth and $L^2>4l_s$ then $\E$ is Bogomolov unstable
by the Bogomolov instability theorem. In this case we also get 
the curve $D$ satisfying assertions of 3.1.(1). This gives better result
for $s$ odd and worse for $s$ even.

\medskip

Now we  state the criterion for jet spannedness at singular points of $X$
indicating necessary changes in the  proof.  
In the theorem $Z=Z_x$ denotes a fundamental cycle of a singularity $(X,x)$
and $\Delta=\Delta_x=f^*K_X-K_Y$ where  $f$ is a minimal resolution of  
the singularity $(X,x)$. Set $\delta _k=-((k+1)Z+\Delta)^2$.

\proclaim{Theorem 3.1'}
Let $L$ be a pseudoeffective Weil divisor on a normal projective surface $X$. 
Assume that $K_X+L$ is  Cartier and generates $(k-1)$-jets but not $k$-jets
at a rational singularity $x$.
\item{(1)} If $L^2>\delta_k$, then there exists a curve 
$D$ containing a Gorenstein cluster $\zeta$ such that $\I_\zeta$ 
contains $m_x^{k+1}$ but it does not contain  $m_x^k$ and such that the map 
$H^0(\O_D(K_X+L))\to H^0(\I _\zeta\O_D(K_X+L))$ is not onto.
In particular, $|\O_D(K_X+L)|$ does not  generate $k$-jets at the point $x$.
Moreover, $L-2D$ is pseudoeffective, numerically nontrivial,
$$LD-\deg \zeta\le 2p_aD-2-K_XD$$
and 
$$LD-{1\over 4}\delta_k \le D^2.$$ 
\item{(2)} If $L^2=\delta_k$, then either there exists a curve $D$
as in (1) or $\epsilon (x, L)=k+1$.
\endproclaim

\pr Instead of the blow up $f$ we  use a minimal resolution of singularity at 
$x$ and then $E$ is replaced by $Z$.
The cluster $\zeta$ is Gorenstein by Theorem 1.5.7, [La3]. The rest of 
the proof is the same.

\medskip
\noindent
{\sl Remarks.}

\noindent
(1) Recall that if $X$ has only Du Val singularities then 
$2p_aD-2\le K_XD+D^2$ (see, e.g., [La3], Corollary 1.3.3; this is just
a simple corollary to the Riemann--Roch theorem as written in 
[La1], Theorem 2.1). Hence in this case the first
inequality in Theorem 3.1', (1) is stronger than $LD-\deg\zeta\le D^2$.
Moreover, in any case this 
inequality bounds discrete invariants of the curve since $p_aD$
and $D(K_X+L)$ are integers.

\noindent 
(2) Note that we have a trivial bound 
$\deg \zeta \le \deg \O_X/m_x^{k+1}=1/2(k+1)(-kZ^2+2)$
(this should be read keeping in mind that $-Z^2=\hbox{emb}\,\dim _xX-1$).
Similarly as in Section 2 one can ask about 
$\max \{\deg\zeta\: \, \zeta{\hbox{  is Gorenstein and }m_x^{n+1}\subset
\I_\zeta} \hbox{ in } \O_{X,x} \}$ but it seems to be quite difficult
question.
\medskip

The following corollary is a Reider type theorem for $k$-jet spannedness
of adjoint divisor for nef $L$
and it considerably improves Corollary 7.5, [Laz] 
(it should be also compared with Proposition 5.7, [Laz];
see Proposition 4.1):

\proclaim{Corollary 3.2}
Let $L$ be a Weil divisor  on a normal projective surface $X$ such
that $K_X+L$ is  Cartier  and $x$ a smooth point of $X$.
If $L$ is nef and $L^2\ge (k+2)^2$ then one of the following holds:
\item{(1)} $\epsilon (x,L)=k+2$ and $L^2=(k+2)^2$,
\item{(2)} $K_X+L$ generates $k$-jets at the point $x$,
\item{(3)} there exists a curve $D$ passing through the point $x$
such that the complete linear system $|\O_D(K_X+L)|$ does not generate
$k$-jets at the point $x$ and such that 
$$LD-l_k\le D^2<1/2LD<l_ k.$$
\endproclaim

\pr
By Theorem 3.1 it is sufficient to prove that if $L$ is nef,
$LD-l_k\le D^2$, $L-2D$ is pseudoeffective and numerically
nontrivial then $LD<2l_k$ and $(L-2D)D>0$.  

We know that $(L-2D)D\ge 0$, i.e., $LD/L^2\le 1/2$, and hence 
by the Hodge index theorem 
$$D^2\le {(LD)^2\over L^2}\le {1\over 2} LD.$$
If we have an equality $D^2={1/2}LD$ then $(L-2D)L=0$ and $L$ and $D$ are 
numerically proportional again by the Hodge index theorem. Using this two 
facts we see that $L-2D$ is numerically trivial, a contradiction.

Therefore $LD-l_k\le D^2<{1/2}LD$, which implies $LD<2l_k$, Q.E.D.

\medskip

In the following theorem $\delta (x,k)$, where $x$ is a point of $X$ and $k$ 
an integer, stands for $(k+1)^2$ if $x$ is smooth and $-(kZ_x+\Delta_x)^2$
if $X$ has a rational singularity at $x$ (and is singular at $x$).

\proclaim{Theorem 3.3}
Let $L$ be a pseudoeffective Weil divisor on a normal projective
surface $X$ such that $K_X+L$ is Cartier and let $k_1$,..., $k_r$ be positive
integers. 
Assume that $X$ has at most rational singularities at distinct points 
$x_1$,..., $x_r$ and the restriction map
$$H^0(\O_X(K_X+L))\to H^0(\O_X(K_X+L)/{m_{x_1}^{k_1}\ot...\ot m_{x_r}^{k_r}})$$
is not onto. If $L^2>\sum _{i=1}^{r}\delta(x_i,k_i)$, 
then there exists a curve 
$D$ passing through $x_1$,..., $x_r$ and such that 
$$H^0(\O_D(K_X+L))\to H^0(\O_D(K_X+L)/{m_{x_1}^{k_1}\ot...\ot m_{x_r}^{k_r}})$$
is not onto. Moreover, $L-2D$ is pseudoeffective, numerically nontrivial, 
$$LD-\sum_{i=1}^r\deg {\O_X/m_{x_i}^{k_i}}\le 2p_aD-2-K_XD$$
and
$$LD-{1\over 4}\sum_{i=1}^r\delta (x_i,k_i)\le D^2.$$ 
\endproclaim

The proof of this theorem is analogous to the proof of Theorems 3.1
and 3.1' and therefore we skip  it. Note that we can also treat the case
$L^2=\sum _{i=1}^{r}\delta(x_i,k_i)$ but the statement of the theorem would be 
more complicated. Similarly as before one can also write down
a more complicated version with
better inequality and some ``bad'' Gorenstein cluster.

Theorem 3.3 can be used
to check $k$-jet ampleness of the adjoint line bundle on surfaces
with at most rational singularities. If the surface is smooth
we have the following theorem (in which we treat also the boundary case):

\proclaim{Theorem 3.4}
Let $L$ be a pseudoeffective divisor on a smooth surface $X$.
Assume that $K_X+L$ is $(k-1)$-jet ample but not $k$-jet ample.
\item{(1)} If $L^2> (k+2)^2$ then there exists a curve $D$ 
such that $\O_D(K_X+L)$ is $(k-1)$-jet ample but not $k$-jet ample, 
$L-2D$ is pseudoeffective, numerically nontrivial and 
$LD-l_k\le D^2$. 
\item{(2)} If $L^2=(k+2)^2$ then either there exists a curve $D$
as in (1) or $k$ is even and there exists a point $x$ such that $K_X+L$ is not 
$k$-jet generated at $x$ and $\epsilon (x,L)=k+2$.
\endproclaim

\pr
This is just a simple generalization of Theorem 3.1
(and the second point of the remark  after this theorem). The only 
thing we need to complete the proof is the inequality
$$\sum_{i=1}^r l_{k_i-1}<l_k$$
if $r\ge 2$ and $\sum _{i=1}^r k_i=k+1$.
This is implicitly proved in the proof of Theorem 2.2. 

\medskip

Similarly as before one can get slightly better theorem if $L$ is nef
(cf.~Corollary 3.2).

\heading 4. Further remarks on jet ampleness \endheading

In this section we assume that $L$ is a divisor on a smooth surface $X$.

The following proposition is a generalization of
Proposition 5.7, [Laz]. The method of the proof is almost the same
but we show it for the convenience of the reader. We think that this
proposition explains appearing of the Seshadri constants in theorems
from Section 3.

\proclaim{Proposition 4.1}
Let $x_1$,..., $x_r$ be distinct points and $k_1$,..., $k_r$
positive integers.  
\item{(1)} If $\sum _{i=1}^r{k_i+1\over \epsilon (x_i,L)}<1$ then 
$H^0(K_X+L)\to H^0(\O(K_X+L)/{m_{x_1}^{k_1}\ot...\ot
m_{x_k}^{k_r}})$ is onto.
\item{(2)} If $\sum _{i=1}^r{k_i+1\over \epsilon (x_i,L)}=1$ 
and $L^2\ge\sum _{i=1}^r(k_i+1)^2$ then 
$H^0(K_X+L)\to H^0(\O(K_X+L)/{m_{x_1}^{k_1}\ot...\ot
m_{x_r}^{k_r}})$ is onto unless $r=1$ and $L^2=(k_1+1)^2$.
\endproclaim

\pr
Let $f\: Y\to X$ be a blow up of $X$ at $x_1$,..., $x_r$
and $E_1$,..., $E_r$ respective exceptional divisors.
It is sufficient to prove that $L'=f^*L-(k_1+1)E_1-...-(k_r+1)E_r$
is nef and big since then $H^1(\O(K_X+L)m_{x_1}^{k_1}\ot ...\ot m_{x_r}^{k_r})
=H^1(K_Y+L')=0$ by the Kawamata--Viehweg vanishing theorem.
Note that
$$L'=\sum_{i=1}^r(f^*L-\epsilon (x_i, L)E_i)+
\left(1-\sum_{i=1}^r {k_i+1\over
\epsilon (x_i,L)}\right)f^*L.$$
Now (1) follows since all the terms on the right are nef and  the last one is
big. (2) is clear if $(L')^2=L^2-\sum_{i=1}^r (k_i+1)^2>0$.
Otherwise $(L')^2=0$ and therefore 
$$(f^*L-\epsilon (x_i, L)E_i)(f^*L-\epsilon (x_j, L)E_j)=0$$
for each pair $i$, $j$. This is impossible if $i\ne j$. Therefore
$r=1$ and $(f^*L-\epsilon (x_1, L)E_1)^2=L^2-(k_1+1)^2=0$, Q.E.D.

\medskip

\proclaim{Corollary 4.2}
Let $A$ be an ample line bundle on a smooth surface $X$. Let
$x_1$, $x_2$,..., $x_r$ denote $r$ distinct points on $X$ and 
$k_1$,..., $k_r$ be some integers such that $k_1+...+k_r=k+1$. 
If $\epsilon (x_i,A)\ge 1$ for $i=1,...,$ $r$, then
$H^0(K_X+nA)\to H^0(\O(K_X+nA)/{m_{x_1}^{k_1}\ot...\ot
m_{x_r}^{k_r}})$ is onto for $n\ge k+2+r$
or for $n\ge k+1+r$ if $A^2>1$.
\endproclaim

\proclaim{Corollary 4.3}
If $A$ is ample, globally generated and $(X,A)\not \simeq (\P^2, \O (1))$ 
then $K_X+nA$ is $k$-jet generated at each point of $X$ for $n\ge k+2$ and
$k$-jet ample for $n\ge 2(k+1)$. 
\endproclaim

\pr
If $A$ is ample and globally generated then $\epsilon (x,A)\ge 1$ for 
any point $x$ of $X$. If $A^2=1$ then the morphism $\varphi$ defined by $A$
is finite and $\deg \varphi \cdot \deg \O _{\varphi (X)} (1)=A^2=1$,
so $\deg\varphi =\deg \O _{\varphi (X)} (1)=1$.
Hence $(X,A)\not \simeq (\P^2, \O (1))$, a contradiction.

Therefore $A^2>1$ and the corollary follows from Corollary 4.2, Q.E.D.

\medskip

\re 
Corollary 4.3 is analogous to Corollary 3.3.(2), [BS2], which says, 
in particular, that if $A$ is very ample 
then $K_X+nA$ is $k$-jet ample for $n\ge k+2$.

\medskip
The following proposition follows from Corollary 3.2: 

\proclaim{Proposition 4.4 }
Let $X$ be a minimal surface of the Kodaira dimension $0$.
Let $A$ be an ample line bundle on $X$ such that $A^2\ge 4$.
If $\epsilon (A)<1$  then there exists an
irreducible curve $D$ such that $AD=1$ and $p_aD=1$. 
Moreover, the curve $D$ has at most double points as its singularities and
if $D$ is singular then $\epsilon (A)=1/2$.
\endproclaim

\re
From the first part of Proposition 4.4 applied for $2A$ it follows that 
$\inf \{ \epsilon (A)\: \hbox { $A$ is ample}\} \ge {1\over 2}$. 
(Now the second part of Proposition 4.4 follows from the inequality
$\epsilon (A)\le AD/mult_xD$ for any $x\in D$.)
It can be also seen by applying usual Reider's theorem since 
$K_X+2A$ is globally generated (note that $A^2\ge 2$ by the Riemann
Roch theorem) and $2\epsilon (A)=\epsilon (K_X+2A)\ge 1$.

This simple remark is related to Problem 3.6, [EL].

\proclaim{Example 4.5} Seshadri constants on K3 surfaces.
\endproclaim

The aim of this example is to give simple proofs of some results obtained in
[BDS]. 

Let $A$ be an ample divisor on a K3 surface $X$.
If $\epsilon (x,A)<1$ for a point $x$ then $L$ is not globally generated.
By  the results of Saint--Donat it follows that $A$ is of the
form $A=aE+\Gamma$, where $E$ is an elliptic curve generating
a free rational pencil $f\: X\to \P ^1$, $\Gamma$ is a $(-2)$-curve with
$E\Gamma=1$, $a\ge 3$. Recall that $\epsilon (y, D)\ge 1$ for
$D$ an ample effective divisor passing through $y$ and smooth at $y$.
It follows that the point $x$ does not lie on $\Gamma$ and
the unique fibre $F$ of f passing through $x$ is singular at $x$.
By definition $\epsilon (x,A)\le AF/\mult _xF=1/\mult_xF$.
Since $\epsilon (x,A)\ge 1/2$ by the remark above, it follows that
$\mult_xF=2$ and $\epsilon (x, A)=1/2$. 

Let $S$ denote a set of singular points of the fibres of $f$.
As a corollary to the above by Proposition 4.1  we get
the following:
\item{1.} 
If $A$ is globally generated then $nA$ is $k$-jet generated for $n\ge k+2$.
\item{2.} 
Otherwise $A$ is of the form $aE+\Gamma$ described above and
$nA$ generates $k$-jets at each point $x\not\in S$  for $n\ge k+2$ and
at each point $x\in S$ for $n\ge 2k+4$.

\medskip

Clearly, $0$-jet spannedness and $0$-jet ampleness are the same thing
but $1$-jet spannedness should not imply $1$-jet ampleness.
However, finding an explicit example of $1$-jet spanned but not
$1$-jet ample line bundle seems to be nontrivial. 
Here we provide such an example:

\proclaim{Example 4.6}
$1$-jet spanned but not $1$-jet ample line bundle.
\endproclaim

Let $X$ be a numerical Campedelli surface with ample $K_X$ and 
$\pi_1^{alg}(X)=\Z_3\op\Z_3$. Theorem 0.3, [La2] together
with Proposition 5.5, [La2] say that there are only 4 degree 2 clusters
which are contracted by $|3K_X|$ and all of them are scheme-theoretical
intersections of unique curves from $|K_X-\tau|$ and $|K_X+\tau|$,
$\tau\in \Tors X-\{0\}$. We will prove that for any surface $X$
all those clusters consist of 2 distinct points.
This can be proved explicitly by using Xiao's construction of such surfaces.
We recall this construction since we need it to do explicit calculations.

Let us choose homogeneous coordinates $([x_0,x_1,x_2],[y_0,y_1,y_2])$
in $\P^2\times\P^2$ and let $\t X_{\lambda}$ be the complete intersection of
two hypersurfaces
$$\sum_{i=0}^{2} x_iy_i=0\quad\hbox{and}\quad
(\sum_{i=0}^{2} x_i^3)(\sum_{i=0}^2y_i^3)-\lambda \prod_{i=0}^2x_iy_i=0.
\leqno(*)$$
For general $\lambda$ this surface is smooth and it is invariant under the 
action of group $G=\Z_3\op\Z_3$ with generators acting by
$$([x_0,x_1,x_2],[y_0,y_1,y_2])\ra([x_1,x_2,x_0],[y_1,y_2,y_0])$$
and
$$([x_0,x_1,x_2],[y_0,y_1,y_2])\ra([x_0,\epsilon x_1,\epsilon^2 x_2],
[y_0,\epsilon^2 y_1,\epsilon y_2]),$$ 
where $\epsilon$ is a primitive cube root of $1$.
The quotient $X_\lambda$ of ${\t X}_\lambda$ by $G$ is a required surface.
Note that $H^0(K_{{\t X}_\lambda})=\bigoplus _{\tau\in \Tors X}H^0(K_
{X_\lambda}+\tau)$ and therefore clusters contracted by $|3K_X|$ are images
of the points defined on ${\t X}_\lambda$ by pairs of equations
$$\sum_{i=0}^2\epsilon^{ai}x_iy_{i+b}=0\quad\hbox{and}\quad
\sum_{i=0}^2\epsilon^{-ai}x_iy_{i-b}=0\leqno(**)$$
for $(a,b)=(0,1),(1,0),(1,2),(1,1)$ 
(the numeration is cyclic modulo 3).

To prove that those clusters consist of distinct points we should
only show that for any $\lambda$ system of equations $(*)$
and $(**)$ have two solutions lying in different orbits of $G$.

We have the following solutions of $(*)$ and $(**)$ lying in different
orbits of $G$:

\item{1.} $([0,0,1],[1,-1,0])$ and $([1,-1,0],[0,0,1])$ for $(a,b)=(1,0)$,
\item{2.} $([0,1,-1],[1,1,1])$ and $([1,1,1],[0,1,-1])$ for $(a,b)=(0,1)$,
\item{3.} $([0,1,-\epsilon^2],[1,\epsilon^2,1])$  and 
$([1,\epsilon,1],[0,1,-\epsilon])$ for $(a,b)=(1,1)$,
\item{4.} $([1,\epsilon^2,1],[0,1,-\epsilon^2])$ and
$([0,1,-\epsilon],[1,\epsilon,1])$ for $(a,b)=(1,2)$.

This shows that the line bundle
$\O_X(3K_X)$ is $1$-jet generated but not $1$-jet ample.

\medskip

By Theorem 0.2, [La2] 
one can expect that for a general surface $X$ in the moduli space of Godeaux 
surfaces with $H_2(X,\Z)=\Z _3$ the line bundle $\O(4K_X)$ is $1$-jet spanned
but not $1$-jet ample. Although an explicit construction 
of such surfaces is known (see [Rd, Section 3]) it is very complicated and
the calculation it involves seems discouragingly large.

\enddocument